\documentclass[11pt]{article}

\usepackage[hyphens]{url}
\usepackage{mathptmx}      
\usepackage{helvet}        
\usepackage{courier}
\usepackage{graphicx}
\usepackage{amssymb}
\usepackage[margin=1in]{geometry}
\usepackage{subfig}       
\usepackage{type1cm}        
\usepackage{makeidx}         
\usepackage{graphicx} 
\usepackage{float}      
\usepackage{multicol} 
\usepackage{subfig}       
\usepackage{type1cm}
\usepackage{tabularx}        
\usepackage{lipsum}
\usepackage{hyperref}   
\usepackage[justification=centering]{caption}   
\newcommand\blfootnote[1]{%
  \begingroup
  \renewcommand\thefootnote{}\footnote{#1}%
  \addtocounter{footnote}{-1}%
  \endgroup}

\begin{document}
\title{Writing Women in Mathematics into Wikipedia}
\author{Marie A. Vitulli \\
University of Oregon\\
Eugene, OR 97403-1222}
\date{October 20, 2017}
\maketitle

\abstract{
In this article I reflect upon the problems connected with writing women in mathematics into
 Wikipedia.  I discuss some of the current projects and efforts aimed at increasing the visibility of women in
  mathematics on Wikipedia.  I present the rules for creating a biography on Wikipedia and relate
   my personal experiences in creating such articles.  I hope to provide the reader with the
    background and resources to start editing existing Wikipedia articles and the confidence to create new articles.  I would also like to encourage existing editors to look out for and protect new articles about women mathematicians and submit new articles.}
\section{Introduction}
The twin problems of the paucity of women subjects and scarcity of women editors\blfootnote{\textcopyright  \hspace{.05 cm} 2017 \hspace{.05cm} Marie A. Vitulli}  on Wikipedia are well known but no 
solutions are on the horizon.  We can nevertheless take small steps to address these problems.   In 
response to my concerns about these problems and in order to gain first-hand experience, I became a Wikipedia editor in 2013.  In this article I will update and expand upon my article \textit{Writing Women in Mathematics
into Wikipedia}, which appeared in the May--June 2014 issue of the \textit{Association for Women in Mathematics Newsletter} \cite{MAV2014}.   Since 2013 some
things have improved on Wikipedia and others have remained more or less the same.  I will first
report on these and then relate some of my personal experiences as a Wikipedia editor. Briefly, this is how this article is organized.  In  Sect. \ref{Problems}, I discuss gender bias on Wikipedia, the underrepresentation of women on Wikipedia, both as editors and subjects of mathematical biographies, and briefly explain my own interest in this subject. 
In  Sect. \ref{CreatingABiography},  I describe the process of creating a biography on 
 Wikipedia and relate my personal 
 experiences of creating pages for prominent women in mathematics on Wikipedia. In  Sect. \ref{Steps}, I  discuss the Wikipedia Year of Science as well as the related projects
\textit{WikiProject Women scientists} and \textit{WikiProject Women in Red}  and 
compare the situation for women in mathematics on Wikipedia in 2014 to the situation in late
 2017.  I also speak 
 about some of the most recent efforts of the Association for Women in Mathematics (AWM) to 
 increase the visibility of women in mathematics in Wikipedia.

\section{Problems for Women on Wikipedia}
\label{Problems}
 
 I have been concerned about the underrepresentation of women in mathematics and the
inadequate recognition of outstanding female mathematicians since the late 1970s when I began
my long career as a university professor of mathematics.  More recently, I have been troubled by
the meager representation of women in mathematics as both subjects of Wikipedia pages and
as Wikipedia contributors, better known on Wikipedia as editors.  

There have been numerous allegations that Wikipedia suffers from systemic gender
 bias with respect to both editors and content \cite{Atlantic2015} and that the attempts to 
 increase women's participation have failed
 \cite{NYT2011},  \cite{HBR2016}.   See the the Wikipedia page \textit{Gender bias on Wikipedia} for additional references and further discussion of that issue \cite{WikiGenderBias}. 
 
 \subsection{The dearth of biographies of women mathematicians on Wikipedia}
One disturbing manifestation of the underrepresentation of women in Wikipedia is the dearth  of 
biographies of notable women. According to \textit{WikiProject Women in Red}, which I will discuss 
in Sect. \ref{Steps}, only 16.36\% of the biographies in English Wikipedia were about women as of August 7, 
2016.  This is up from just over 15\% in November 2014. Wikipedia has guidelines on academic
 notability (also called the professor test) that a subject must meet to merit a Wikipedia page. The
  first two of the nine criteria of academic notability are: the personÕs research has made significant
   impact in their scholarly discipline, broadly construed, as demonstrated by independent reliable
    sources; and the person has received a highly prestigious academic award or honor at a
     national or international level. Some of the criteria are subjective in nature and it is not surprising that some of the American mathematicians featured on Wikipedia are more notable than
      others.
      
In 2013 the Sociology, Science, Tech \& Society Program of the National Science Foundation awarded a grant to Yale University [1] to explore potential gender differences in indicators of academic notability, the existence of networks that help contributors mobilize resources for content creation, gate keeping processes that result in challenges to content associated with women and scholarship, and the unintended consequences of WikipediaÕs policies themselves.

When I first perused the page for the \textit{Category:Women mathematicians} (formerly called 
\textit{Category: Female mathematicians}) several years ago I was disheartened to see how many 
prominent women lacked a Wikipedia page. I downloaded a copy of the category in August of
 2015 and again in August of 2016. I was pleased to discover that between 2015 and 2016 there
  were 73 new biographies of women mathematicians on Wikipedia. I also keep an eye on the
   women included in the \textit{Category:American mathematicians}. As of August 30, 2016 there 
   were
    37 pages (i.e., people) in the \textit{Category:American women mathematicians}. If we assume every
     woman in the \textit{Category:American women mathematicians} is also in the \textit{Category:American 
     mathematicians} (which is not the case), the percentage of women in the \textit{Category:American
      mathematicians} was about 17\%, at that time. As of October 12, 2017 there were 73 women listed in the \textit{Category:American women mathematicians}; the number of women listed nearly doubled between August 2016 and October 2017, in part, due to the efforts I discuss in Sect. \ref{Steps}.  There is also the \textit{Category:Women statisticians}, which contained 82 pages as of October 12, 2017.  According to American Mathematical Society (AMS) data from the Annual Survey, we know that 29\% of the U.S. citizens to earn PhDs from U.S. mathematics departments between 1991 and 2015 were women. As per the 2015 AMS Annual Survey, about 23\% of all doctoral mathematics faculty were women and 19\% of the Public or Private Large doctoral mathematics faculty were women.

\subsection{Writing Wikipedia pages that will survive}  
A serious impediment to having more biographies of women in mathematics on 
Wikipedia is that though anyone can become a Wikipedia editor and post an article, not all of the articles that are created 
will ``survive'' on Wikipedia.          There are various ways that an article 
         can be deleted.  If no opposition to the deletion is expected, \textit{any} editor 
         can tag the article for \textit{Proposed deletion} (\textbf{PROD}), which lists the article in the 
         \textit{Category:Proposed deletion}. Any editor (including the creator of the article) may 
         object to the deletion by simply removing the tag and giving an explanation. My own experience is that the author of an article should be cautious about removing the \textbf{PROD} tag without giving a well-reasoned explanation that is backed up by citations.  An article 
         marked for proposed deletion is checked by an Administrator, who may remove the \textbf{PROD} tag or delete the article.  An article may also be nominated for the page 
         \textit{Articles for deletion} \textbf{(AfD}), which is where Wikipedians discuss whether or not 
         an article should be deleted.\footnote{A Wikipedian is a volunteer who writes or edits Wikipedia articles.}  An article listed here is usually discussed for at least seven 
         days, at the end of which the deletion process proceeds on the basis of community 
         consensus.  I will discuss my  personal experience with the \textbf{PROD} tag in Sect. \ref{CreatingABiography}.
         
         Several women who have created biographies of women in mathematics on Wikipedia recommend that the author of a new page have some editors on standby to affirm the article created by making minor edits soon after the article goes live.  Edits can be made by less senior editors without logging in to their Wikipedia accounts; only the IP address of the computer they were using will appear on the History page of the article.  It is even better if a page creator has  an active Administrator looking out for a newly created  article.
         
         There is a hierarchy of ``functional capacity'' 
        among the Wikipedia participants \cite{FormalOrg}.  Wikipedia co-founder Jimmy Wales is 
        the ultimate arbiter on Wikipedia, however, he has usually deferred his authority to the 
        approximately 34 Bureaucrats or \textit{Crats}, the approximately 700 active Administrators or
         \textit{Admins}, and another group called the Arbitration Committee or \textit{Arbcom}. The membership of the
         Arbitration Committee
         changes annually and currently consists of about 14 Wikipedians, all of whom were formerly 
         Administrators.  The rights and privileges of the \textit{Crats}, \textit{Admins} and \textit{Arbcom} are 
         described on the Wikipedia pages for these groups \cite{WikiBureaucrats}, \cite{WikiAdmins}, \cite{WikiArbCom}.

  \subsection{The scarcity of women editors on Wikipedia}
  Women are also substantially underrepresented on Wikipedia as editors.   A survey \cite{UNUMerit2010} conducted by 
   the United Nations University and Maastricht University in collaboration with the
    Wikimedia Foundation demonstrated that in 2008 only 12.64\% of Wikipedia editors worldwide 
    were women.   A follow up survey in 2011 showed that things hadn't improved:  9\% of editors 
    worldwide and 15\% of  those from the  US  were women.   
    
    A few years after the initial study, the 
    Wikimedia Foundation announced the goal of raising the proportion of female editors  worldwide  to 
    25\% by 2015. In 2014 Wikimedia Foundation founder Jimmy Wales announced that Wikipedia 
    had completely failed to
         reach this goal \cite{BBCNews2014} despite launching several initiatives.  This is particularly disturbing in light of the fact that the percentage of adult Americans who use
         Wikipedia for information increased from 25\% in February 2007 to 42\% in May 2010, 
         according to a study conducted by the Pew Research Center \cite{Pew2011}.   The Pew 
         study also showed that about 50\% of the adult female Internet users and 56\% of the male 
         adult Internet users look to Wikipedia for information.  
         
         In Sect. \ref{Steps}  I will discuss a recent partnership between the AWM and the Wiki Education Foundation to help lessen both the scarcity of women mathematicians as subjects and editors on Wikipedia.

\section{Creating a biography on Wikipedia}
\label{CreatingABiography}
\subsection{The basics: setting up an account, academic notability, and using reliable sources} 

Every editor is encouraged to sign up for an account on Wikipedia's home page \cite{WikiHome}, although this is not mandatory.  An editor should
 care in choosing a username because this name  will be displayed on the editor's
 User page and on the History page of every article created or edited. \footnote{My username is Mvitulli.}  
 After  logging in to a Wikipedia account the User will see a link to her Sandbox at the top of the page.  
 This is where a User can create a draft of an article.  If an editor wants to work on two articles at the 
 same time, she can create a second Sandbox by creating a page entitled  
 User:YourUsername/sandbox2.  \footnote{For example, if one searches for the page User:Mvitulli/sandbox2,
 a template for creating a Wikipedia page on a woman in mathematics will be found.} In this way, an editor can 
 create as many Sandboxes as  desired and with whatever names are desired.

A biographer must remain neutral and there can't be any conflict of interest with the subject.  A current Wikimedia administrator recommends that a biographer start her article with one sentence describing what her subject accomplished and another sentence pointing out why her subject and the subject's work are important. The subject must be notable and every claim made about her must be backed up by a credible citation. Wikipedia has criteria for notability, in general, and for academic notability, in particular.  Wikipedia's guidelines on academic notability (also known as the \textit{professor test}) stipulate that a professor is notable provided that she meets at least one of nine criteria listed on the Wiki page for academic notability \cite{notability_academic}.  Here are the first few criteria from that list.
\begin{enumerate}
\item The person's research has made significant impact in their scholarly discipline, broadly construed, as demonstrated by independent reliable sources.
\item The person has received a highly prestigious academic award or honor at a national or international level.
\item The person is or has been an elected member of a highly selective and prestigious scholarly society or association (e.g., a National Academy of Sciences or the Royal Society) or a fellow of a major scholarly society which reserves fellow status as a highly selective honor (e.g., Fellow of the IEEE)
\end{enumerate}

There are additional specific criteria notes on the Wikipedia page on academic notability.  A 
contributor must demonstrate notability using \textit{reliable sources}. In general, reliable sources 
include mainstream news media and major academic journals, both of which are subject to some 
sort of editorial control. For Wikipedia purposes,  published original research articles are 
considered primary sources; expository articles or summaries of the current state of affairs in a 
particular research area are considered to be secondary sources.  Unlike the case in academic 
publishing, Wikipedia administrators frown on the use of too many primary sources, as I found out the hard way. Wikipedia also has a page on Biographies in their Manual 
of Style  and a special page on Biographies of living persons  \cite{WikiBiographies}, \cite{WikiManualofStyle}.

\subsection{Finding women in mathematics on Wikipedia}
There are many ways to find lists of women mathematicians on Wikipedia. Wikipedia includes 
both an article entitled \textit{List of women in mathematics} and  lists of women with pages 
that have been designated as belonging to the \textit{Category:Women mathematicians} or to the 
\textit{Category:American women mathematicians}. When an editor inserts  the tag 
[[Category:Women mathematicians]] at the end of her article, the name the subject is  automatically added to the page for that category. There are many category tags 
 that may be used; look at a few articles to see which may be appropriate.


\subsection{My experiences writing biographies on Wikipedia}

I attended my first edit-a-thon on writing women into Wikipedia at the University of Oregon in March of 2013. This event was part of the Wiki Women's History Month events.   The edit-a-thon  was led by Sarah Stierch, who was then a
 program evaluation coordinator for the Wikimedia Foundation and a former Wikipedian in Residence at the Archives of American Art and the Smithsonian Institution Archives.

During the edit-a-thon I decided to write an article about Susan Montgomery.  Montgomery was 
the 2011 AWM Noether Lecturer and in 2012 she was selected as both an American Association for the Advancement of Science (AAAS) Fellow and an  
AMS Fellow in the inaugural class.  I felt she clearly merited a page on Wikipedia.  I began 
working on her page during the edit-a-thon and finished up the article a couple of weeks later.  My 
first experience writing for Wikipedia was extremely frustrating, but in the end I am glad that I undertook this project.  I learned that there are many Wikipedia conventions that a contributor 
must follow.  I am still learning about  these conventions and the culture of Wikipedia.

The first thing I did during the edit-a-thon was to edit the existing Wikipedia page on the Noether
 Lecture  so that the lecturers for 2011 -- 2013, including Montgomery, were listed on that page. I inserted square brackets around Montgomery's name so she appeared as a red link.  A biography has a better chance of surviving on Wikipedia if the subject is already mentioned on
  an existing Wikipedia page.  I followed all the above-mentioned guidelines and asked Stierch to
   read my article before I published it.  In spite of all my precautions my article was proposed for
   deletion with the \textbf{PROD} tag within ten minutes after it went live by a novice editor who made this remark.

\begin{quotation}
It is proposed that this article be deleted because of the 
following concern: This looks like a case of shameless (self) 
promotion. Hey, look at me: I teach math and I want my CV on 
Wikipedia. 
\end{quotation}

Stierch responded to the novice editor by saying the biography was in fact a "good faith article by a new editor (not the subject)" and removed the \textbf{PROD} tag.  The article survived only because of Stierch's intervention. \footnote{Stierch, whose username today is Missvain, is still an active Wikipedia \textit{Admin}. 
Regrettably, Stierch lost her Wikimedia Foundation job in late 2013 when she was discovered editing for pay \cite{StierchFired}.   Wikipedia has a conflict of interest policy on paid editing and a page on this policy \cite{Conflict}.  In general there is a prohibition on paid editing but it isn't uniformly enforced. There is another page \cite{Paid Editing} discussing this.  I have heard that the Ivy League schools have public relations staff who create pages for distinguished faculty; the staff are supposed to acknowledge their employer in their Talk pages. I came across one such acknowledgment on the User page of Dominic McDevitt-Parks (User:Dominic) who is a paid employee of the US National Archives and Records Administration and edits on Wikipedia for them. I suspect that there are many other paid editors who make no such acknowledgement.}

At Montgomery's suggestion I started another Wikipedia biography in May of 2014, this one on Georgia Benkart.  Benkart was AWM President during 2009-2011, was elected to the inaugural class of Fellows of the AMS in 2013, was selected to deliver the  AWM-AMS Noether Lecture in 2014 and the Emmy Noether Lecture at the International Congress of Mathematicians  in Seoul, Korea in 2014.   She has published over 100 journal articles and has co-authored three \textit{Memoirs of the American Mathematical Society}.  Benkart  clearly met the Wikipedia criteria for academic notability.   I created the page in my Sandbox and had  Stierch both read  it and later monitor it when the article was published it on Wikipedia.  I also asked Emily Temple-Wood to watch the page and make sure it wasn't deleted.

Soon after the Benkart biography went live, it was tagged for having too many primary sources and 
not enough secondary sources. Wikipedia recommends that a biography  should have at least as many secondary sources in the References section  as  publications of the subject, which are regarded as primary sources.  I addressed that criticism and removed the tag after thanking the 
editor involved  for his comments and explaining that most of the sources were now secondary.   
The article was then tagged for using \textit{weasel} words: vague phrasing that often accompanies 
biased or unverifiable information \cite{WikiWeasel}. The Wikipedia administrator asserted the following. 

\begin{quotation}
The whole tone of the article is written putting her on a pedestal.
\end{quotation}
After asking Emily Temple-Wood and 
Sarah Stierch about weasel words and looking at the Talk page for the article, I found out that the administrator was objecting to my claims that Benkart was a ``distinguished" mathematician who is an ``international leader" and a 
``renowned" teacher, and for asserting that one of  Benkart's joint papers  became ``one of the building blocks" of the classification of toroidal rank-one Lie algebras.  Since I am not an expert in Benkart's research 
  area (the structure and representation theory of Lie algebras), I asked  Montgomery for 
  help.  Montgomery enlisted the help of Efim Zelmanov, who wrote up the paragraphs that 
  described the ``importance'' of Benkart's work (I couldn't help slipping in the weasel word ``importance'').  In this case an 
  administrator who has virtually no background in mathematics objected to a
  Fields medalist's description of Benkart's work.  As Stierch pointed out to me, this is the reality of 
  the wacky world of Wikipedia.  If an editor can't prove her assertion with a reliable secondary source, 
  then  any other editor can remove the relevant sentence without question.  Stierch removed the 
  uncited weasel words from my article and the tag was eventually removed. As of this writing, the 
  article has survived on Wikipedia.
  
 Other mathematicians have also had a difficult time when writing biographies of women 
 mathematicians.   For example, Sormani has created several biographies, originally with the username ``Sormani."  Beginning around 2007, the articles she 
 wrote were deleted by 
 other editors and \textit{Admins}.  She then started writing articles anonymously without logging 
 on to the Wikipedia site and had more success.  We do not have any data, numerical or anecdotal, on the difficulties editors have in creating or editing pages on male mathematicians.

\section{Steps towards a solution}
\label{Steps}

During the past several years there have been a number of successful attempts to increase the visibility of women mathematicians on Wikipedia. In 2007 the Wikipedia page \textit{Noether Lecture} was created by AWM member Christina Sormani. The Noether Lecture is lecture series (and award) that honors women ``who have made fundamental and sustained contributions to the mathematical sciences.'' The AWM established the lecture in 1980 and in 2013 it was renamed the AWM-AMS Noether Lecture. The Noether lecturers were all listed on the page \textit{Noether Lecture}, some as red links, that is, people (or topics) mentioned on Wikipedia but for whom Wikipedia pages do not exist\footnote{A Wikipedia editor can easily create a red link on Wikipedia by 
 surrounding [[name]] or [[topic]]  with double square brackets as indicated.}.  I am happy to report that as of today all of the Noether lecturers listed on the page Noether Lecture appear as blue links, that is, people or topics with pages on Wikipedia \cite{WikiNoetherLectures}.


\subsection{The WikiProjects Women scientists and Women in Red}
The \textit{WikiProject Women scientists} \cite{WP Women scientists} began in 2012 and is dedicated to ensuring the quality
and quantity of biographies of women in science on Wikipedia.  Emily Temple-Wood (who edits on
 Wikipedia as Keilana) is a co-founder of this project; she started  writing Wikipedia pages when 
 she
 was 12 and continues to write pages and confront gender bias on Wikipedia, although her pace has slowed a bit 
 since she  began 
 medical school at Midwestern University in the fall of 2016.  Temple-Wood has written over 100 pages on women scientists.  She became an \textit{Admin} in 2007 and is currently part of the Arbitration Committee.  There is a plethora of information on the project page \cite{WP Women scientists} for \textit{WikiProject Women scientists},  including suggestions of how an editor  can help the project. There are links to pages that are helpful to Wikipedia editors and a list of more than 80 members of the project. \footnote{The list of members begins with the usernames of Emily Temple-Wood and Sarah Stierch (Keilana and SarahStierch, respectively); I am member 84 near the bottom of the list.}   
 
The objective of the \textit{WikiProject Women in Red} 
 \cite{WP Women in Red} is to
 turn red links on Wikipedia (names and topics that appear on Wikipedia but do not have Wikipedia pages) into blue links (names or topics
 with Wikipedia pages). According to \textit{WikiProject Women in Red},  only 16.36\% of the  biographies in English Wikipedia
were   about women as of August 7, 2016.  This is up from just over 15\% in November
 2014.  Rosie Stephenson-Goodknight (who edits on Wikipedia as Rosiestep) is a co-founder of
 this WikiProject.   \textit{WikiProject Women in Red} hosts edit-a-thons and socializes their scope and objective  
 via social media. \textit{WikiProject Women in Red} has a page on Mathematics \cite{WomenInRedMath} that is organized by country and
 lists missing articles on women who are (or have been) notable for their contribution to 
 mathematics in academics, business, economics, politics, research, government or the social 
 sector.  I made many additions to the list and encourage others to do so as well.
 
 Temple-Wood and Stephenson-Goodknight were named as 2016 co-Wikipedians
 of the Year by Wikipedia Trustee Jimmy Wales in recognition of their work on harassment on
 Wikipedia along with their collective efforts to expand coverage of notable women on Wikipedia.

\subsection{Other efforts}
 In the recent past, pages for women 
elected as Fellows of the American Mathematical Society (AMS) or fellows of other major professional organizations have been 
routinely created. 
Formerly, on the Talk page of the \textit{WikiProject:Women scientists}, in the topmost section entitled ``Some missing mathematicians", there was a list of 
women who had been elected Fellows of the AMS but didn't have pages;  soon 
after the list appeared pages were created for most of these women by User:Brirush and User:David Eppstein.  After checking this Talk page recently, I discovered that there were five women who were elected Fellows of the AMS in 2017 Class who did not have pages; I  listed their names on the Talk page of \textit{WikiProject:Women scientists}, on the home page for \textit{WikiProject Women in Red/Mathematics}, and on the home page for \textit{Wikipedia:Year of Science/Computer science, technology and math}.  I am happy to report that as of October 12, 2017 all of these prominent women have Wikipedia pages.  This wasn\rq t done for the class of 2018 AMS Fellows.  I created a page for 2018 Fellow Rachel Justine Pries. Almost immediately, another Wikipedia editor, User:FormalDude, attached the PROD tag to the article and proposed it for deletion. I successfully defended the notability of Pries and had a team of editors watch the page; to date the page has survived. Pages for 2018 Fellows Antonella Grassi and Hema Srinivasan were created by Ursula Whitcher and David Eppstein, respectively.  At this point all the women in the class of 2018 AMS Fellows had Wikipedia pages.  There are still women who were elected SIAM Fellows that do not have pages.

In the summer of 2016 Jami 
 Mathewson of the Wiki Education Foundation created the \textit{Category:Awards and prizes of 
 the Association for Women in Mathematics} \cite{WikiCategoryAWMPrizesAwards}.  As of October 12, 2017 there were four pages in this category:
 AWM/MAA Falconer Lecturer, Louise Hay Award, M. Gweneth Humphreys Award, and Noether 
 Lecture.  
 Some of the mathematicians listed on the first three pages appear as red links.   The Wikipedia page for 
 the Association for Women in Mathematics lists all past and current AWM presidents.  Despite the 
 fact that all AWM presidents have Wikipedia pages, and hence appear as blue links, Mathewson reports that only one page (the 
 page for  Ruth Charney) has been rated as ``B-class" by the WikiProjects Biography and Women 
 scientists, designating it as one of Wikipedia's higher-quality 
 articles.\footnote{ To  see these ratings,  click  on the tabs for the Talk pages of the articles.} \footnote{We  discuss the 
 WikiProject Women scientists\cite{WP Women scientists} in more
  detail in Sect. \ref{Steps}} Many of the pages  for former AWM presidents   are \textit{stubs}, that is, pages that exist but need to be expanded to meet Wikipedia standards. 
  To be fair, many of the pages are not rated. My point in mentioning this assessment is that there
   is still a lot of work to do for fledgling editors who may not yet want to create a page from scratch. 

\subsection{The year of science}   
The year  2016 was the Wikipedia Year of Science \cite{WikiYearScience}, an unparalleled initiative to improve scientific
literacy through Wikipedia.  A fundamental part of the Year of Science was the Women in Science
initiative, whose stated goal was to create or improve coverage of women scientists on Wikipedia.
Three groups, \textit{WikiProject Women scientists} \cite{WP Women scientists},  
\textit{WikiProject Women in Red} \cite{WP Women in Red} and the Wiki Education Foundation \cite{WikiEducationFoundation},  
 joined forces to run virtual edit-a-thons during the year to create and improve articles about 
women
scientists in diverse fields, emphasizing a different area each month between February and
December.  Each month had a theme with its own page that was designed to encourage 
Wikipedians to collaborate, share successes, identify articles that need to be written, find 
subject-specific
sources, or initiate any activity that improves coverage of or celebrates women in science.  
During the first several
months  of 2016, the fields of zoology, chemistry, environmental science, psychology, astronomy, 
plant biology, sociology, and medicine were highlighted, whereas the year ended with the spotlight 
on linguistics and physics \cite{WikiYearScience}.  The month of
October  highlighted  computer science, technology, and mathematics and the corresponding Wikipedia page \cite{OctWomenSci}
included a section entitled \textit{Biographies of mathematicians to create or improve}, which I  
authored.
In this section one can find a list of prominent women in mathematics who either don't have
Wikipedia pages or whose pages are stubs.  The page also contains a section labeled
\textit{Association for Women in Mathematics} that lists AWM awards, prizes, and lectures that 
either don't have pages or some of whose recipients don't have pages.  This page is a good starting point for someone who wants to edit or create a page enhancing the visibility of women in mathematics on Wikipedia.  

\subsection{AWM's partnership with the Wiki Education Foundation}

While writing this article I learned of an unprecedented partnership between the AWM and the 
Wiki Education Foundation.   There is considerable need for such a partnership, particularly since there don't 
seem to be many mathematicians involved in the \textit{WikiProject Women scientists} and the \textit{WikiProject Women in 
Red}. As I said earlier, I added my name as a member of the \textit{WikiProject Women scientists} but I didn't see any other usernames that I recognized.  In the past, I 
have sought out Stierch and  Temple-Wood for advice and for help protecting the pages I 
wrote.  Neither of these allies has much time to devote to Wikipedia today.   In response to a query I posted on Temple-Wood's Talk page, Susan Barnum, a Public Services Librarian who lives is El Paso, TX  and has the username Megalibrarygirl, offered help with copyediting and references.  It is crucial that we 
add women editors who can create pages for prominent women in mathematics.  But this effort will 
fail unless there is a network of Wikipedia editors and administrators who watch the pages that we 
create.  Perhaps the partnership can help build such a network.  

AWM  recently appointed two AWM-Wikipedia Visiting Scholars: Tarini
            Hardikar, a senior mathematics-chemistry double major at Colby College, and Sara Del
             Valle, a Ph.D. applied mathematician who is currently a member of the Los Alamos National Laboratory epidemic modeling research team.    Hardikar and Del Valle will be working with Ryan McGrady of the Wiki Education Foundation.  
             Visiting Scholars  are 
           experienced Wikipedians who are capable of improving the quality of articles to B-class or better.  Hopefully these two efforts will increase the numbers of women in mathematics who are
            Wikipedians and increase the quantity and quality of articles on women in mathematics on
             Wikipedia. 

In another unprecedented move, Ursula Whitcher organized a  
Wikipedia women in mathematics edit-a-thon at the 2017 AWM Research Symposium at University of California, Los 
Angeles. The edit-a-thon was supported by a grant from the Wiki Education Foundation.   Mathewson  and I assisted Whitcher during the edit-athon.  This event and ones like it were are a wonderful resource for novice editors.  Ten new pages on women mathematicians were created during the edit-a-thon, including a page I created on 2017 AMS Fellow Mei-Chu Chang.  Please view and help improve these pages by visiting the Wikipedia page for the edit-a-thon \cite{WikiAWMSymposium}. Whitcher will lead a AWM Wikipedia Meetup at the 2018 Joint Mathematics Meetings in San Diego, CA.
 
\section{Summary}                      
There are many more Wikipedia articles on women mathematicians today than when I first became a Wikipedia editor in early 2013. However, many of these articles are stubs and not a single women mathematician appears in the 10 articles allotted to mathematicians in \textit{Wikipedia:Vital Articles}. WikiProject Women scientists is attempting to have the page on Emmy Noether included in this prestigious list. Although all of the women who were elected AMS Fellows now have Wikipedia pages, several of the women who were elected SIAM Fellows or who have received awards from national organizations still do not have pages. The percentage of women editors on Wikipedia remains dismally low
I have edited many biographies of women mathematicians and will continue to do so.

Even though writing a new biography is at times frustrating and is always time-consuming, I am 
glad that I created  biographies of four prominent women mathematicians and hope to write more.  
Anyone who wants to become 
a Wikipedia editor can find a useful list of Wikipedia help pages on the bottom of the 
Project page for the \textit{WikiProject Women scientists} \cite{WP Women scientists}. Placing \{WPWS\}
at the top of the Talk page of an article will result in  the article being automatically added to \textit{Category:WikiProject women 
scientists} and the other members of the WikiProject will be able to keep an eye on the page.  I 
hope you will decide to take the plunge. The entire community of women in mathematics will be 
looking forward to reading new articles on notable women in mathematics.

\section{Acknowledgements}
I would like to acknowledge and applaud the efforts of Emily Temple-Wood and Rosie Stephenson-Goodknight for their work on harassment on Wikipedia along with their collective efforts to expand coverage of notable women on Wikipedia. On behalf of the greater mathematics community I offer them both a great big virtual thank you.

\end{document}